\documentclass[12pt]{article}

\usepackage[ruled,vlined,linesnumbered]{algorithm2e}
\usepackage{amsfonts,amsmath,tikz}

\usepackage{graphicx}
\usepackage{amsmath,amsfonts,amssymb,latexsym}
\usepackage{enumerate}
\usepackage{setspace}
\usepackage{color}
\usepackage[cp1250]{inputenc}
\usepackage[width=145mm,height=210mm,left=35mm,foot=10mm]{geometry}

\newtheorem{thm}{Theorem}[section]

\newtheorem{conj}[thm]{Conjecture}
\newtheorem{cor}[thm]{Corollary}

\newtheorem{lemma}[thm]{Lemma}

\newcommand{\qed}{\hfill ~$\square$\bigskip}
\newcommand{\proof}{\noindent{\bf Proof.} }

\newcommand{\cp}{\,\square\,}
\tikzstyle{vertex}=[circle, draw, inner sep=0pt, minimum size=6pt]

\newcommand{\ggrz}{\gamma_{gr}^{\textrm{Z}}}
\newcommand{\ggr}{\gamma_{gr}}
\newcommand{\ggrt}{\gamma_{gr}^{t}}
\newcommand{\ggrl}{\gamma_{gr}^{L}}

\begin{document}

\title{On Grundy total domination number in product graphs}
 
\author{
Bo\v stjan Bre\v sar $^{a,b}$ \and Csilla Bujt{\'a}s $^{c}$ \and
Tanja Gologranc $^{a,b}$ \and Sandi Klav\v zar $^{d,a,b}$ \and Ga\v sper
Ko\v smrlj $^{b,f}$ \and Tilen Marc $^{b,d}$ \and Bal{\'a}zs Patk{\'o}s $^{e}$ \and Zsolt Tuza $^{c,e}$ \and M{\'a}t{\'e} Vizer $^{e}$ }

\maketitle

\begin{center}
$^a$ Faculty of Natural Sciences and Mathematics, University of Maribor, Slovenia\\
\medskip

$^b$ Institute of Mathematics, Physics and Mechanics, Ljubljana, Slovenia\\
\medskip

$^c$ Faculty of Information Technology, University of Pannonia, Veszpr\'em, Hungary\\
\medskip

$^d$ Faculty of Mathematics and Physics, University of Ljubljana, Slovenia\\
\medskip

$^e$ Alfr\'ed R\'enyi Institute of Mathematics, Hungarian Academy of
Sciences, Budapest, Hungary
\medskip

$^f$ Abelium R\&D, Ljubljana, Slovenia\\
\medskip

\small \texttt{bostjan.bresar@um.si, \{bujtas,tuza\}@dcs.uni-pannon.hu, sandi.klavzar@fmf.uni-lj.si, tilen.marc@imfm.si, patkos@renyi.hu, gasperk@abelium.eu, \{tanja.gologranc,vizermate\}@gmail.com}
\medskip

\end{center}
\maketitle

\begin{abstract}
A longest sequence $(v_1,\ldots,v_k)$ of vertices of a graph $G$ is a Grundy total dominating sequence of $G$ if for all $i$, $N(v_i) \setminus \bigcup_{j=1}^{i-1}N(v_j)\not=\emptyset$. The length $k$ of the sequence is called the Grundy total domination number of $G$ and denoted $\gamma_{gr}^{t}(G)$. In this paper, the Grundy total domination number is studied on four standard graph products. For the direct product we show that $\gamma_{gr}^t(G\times H) \geq \gamma_{gr}^t(G)\gamma_{gr}^t(H)$, conjecture that the equality always holds, and prove the conjecture in several special cases. For the lexicographic product we express $\ggrt(G\circ H)$ in terms of related invariant of the factors and find some explicit formulas for it. For the strong product, lower bounds on $\ggrt(G \boxtimes H)$ are proved as well as upper bounds for products of paths and cycles. For the Cartesian product we prove lower and upper bounds on the Grundy total domination number when factors are paths or cycles. 
\end{abstract}


\noindent
{\bf AMS Subj.\ Class.\ (2010)}: 05C69, 05C76

\sloppy
\section{Introduction}

In~\cite{bgmrr-2014} the Grundy domination number of a graph was introduced as the largest value obtained by a greedy domination procedure. When domination is replaced by total domination (i.e., instead of closed neighborhoods, one considers open ones) the Grundy total domination number arises. The latter was introduced recently in~\cite{bhr-2016}, and was further studied in~\cite{bknt2017}.
   
Let $S=(v_1,\ldots,v_k)$ be a sequence of distinct vertices of $G$. The corresponding set $\{v_1,\ldots,v_k\}$ of vertices from the sequence $S$ will be denoted by $\widehat{S}$. The initial segment $(v_1,\dots,v_i)$ of $S$ will be denoted by $S_i$. A sequence $S=(v_1,\ldots,v_k)$, where $v_i\in V(G)$, is a {\em (legal) open neighborhood sequence} if, for each $i$,
\begin{equation}
\label{eq:defGrundyT}
N(v_i) \setminus \bigcup_{j=1}^{i-1}N(v_j)\not=\emptyset\,.
\end{equation}
Note that if $G$ is without isolated vertices and $S$ is a maximal open neighborhood sequence of $G$, then $\widehat{S}$ is a total dominating set of $G$ (initially, total dominating sequences were introduced just for graphs with no isolated vertices). 

If~\eqref{eq:defGrundyT} holds, each $v_i$ is said to be a {\em legal choice for Grundy total domination}. We will say that $v_i$ {\em totally footprints} the vertices from $N(v_i) \setminus \bigcup_{j=1}^{i-1}N(v_j)$, and that $v_i$ is the {\em total footprinter} of any $u\in N(v_i) \setminus \bigcup_{j=1}^{i-1}N(v_j)$. Any legal open neigborhood sequence is called a {\em total dominating sequence}. For a  total dominating sequence $S$ any vertex in $V(G)$ has a unique total footprinter in $\widehat{S}$. An open neighborhood sequence $S$ in $G$ of maximum length is called a {\em Grundy total dominating sequence} or {\em $\ggrt$-sequence}, and the corresponding invariant the {\em Grundy total domination number} of $G$, denoted $\ggrt(G)$. 

If condition~\eqref{eq:defGrundyT} is replaced with 
\begin{equation}
\label{eq:defGrundy}
N[v_i] \setminus \bigcup_{j=1}^{i-1}N[v_j]\not=\emptyset\,,
\end{equation}
then one speaks of a {\em legal choice for Grundy domination}. All the definitions from the above paragraph can now be restated just by omitting ``total'' everywhere. In particular, the maximum length of a legal dominating sequence in $G$ is the {\em Grundy domination number} of $G$ and denote it by $\gamma_{gr}(G)$. This invariant was studied for the first time in~\cite{bgmrr-2014}.

In \cite{bbgkkptv2016} the Grundy domination number of grid-like, cylindrical and toroidal graphs was studied. More precisely, the four standard graph products of paths and/or cycles were considered, and exact formulas for the Grundy domination numbers were obtained for most of the products with two path/cycle factors.   
In this paper we follow this work and investigate the Grundy total domination number on the standard graph products. The main difference in proving the results for the Grundy total domination number comparing to related results from~\cite{bbgkkptv2016} is that the analogue of \cite[Lemma 1]{bbgkkptv2016} does not hold for Grundy total dominating sequences, therefore other techniques are required. In particular, as in the case of total domination (cf.\ the monograph~\cite{heye-2013}), also in the Grundy total domination natural connections with hypergraphs are applicable.

In \cite{bbgkkptv2017} a strong connection between the Grundy domination number and the zero forcing number of a graph was established. (The zero forcing number is in turn very useful in determining the minimum rank of a graph~\cite{AIM}.) Lin~\cite{L-2017} noticed a similar connection between the Grundy total domination number and the skew zero forcing number of graphs, where the skew zero forcing number was introduced in~\cite{IMA2010} and is denoted by $Z_{-}(G)$. (For some recent related results see~\cite{Aazami2008,ajps-2016,BBF2013,BFF-2016+,D2014,K2015,masu-2016}.) More precisely, Lin proved that if $G$ is a graph, then $\ggrt(G)=|V(G)|-Z_{-}(G)$. As a consequence, $\ggrt(G)\le $mr$_0(G)$, where mr$_0$ is a variation of the minimum rank. 

The vertex set of each of the four standard graph products of graphs $G$ and $H$ is equal to $V(G)\times V(H)$. In the \emph{direct product} $G\times H$, vertices $(g_{1},h_{1})$ and $(g_{2},h_{2})$ are adjacent when $g_1g_2\in E(G)$ and $h_1h_2\in E(H)$. In the \emph{lexicographic product} $G\circ H$ (also denoted in the literature by $G[H]$), vertices $(g_{1},h_{1})$ and $(g_{2},h_{2})$ are adjacent if either $g_{1}g_{2}\in E(G)$, or $g_{1}=g_{2}$ and $h_{1}h_{2}\in E(H)$. In the \emph{strong product} $G\boxtimes H$, vertices $(g_{1},h_{1})$ and $(g_{2},h_{2})$ are adjacent whenever either
$g_1g_2\in E(G)$ and $h_1=h_2$, or $g_1=g_2$ and $h_1h_2\in E(H)$, or $g_1g_2\in E(G)$ and $h_1h_2\in E(H)$.
Finally, in the \emph{Cartesian product} $G\cp H$, vertices $(g_{1},h_{1})$ and $(g_{2},h_{2})$ are adjacent if either
$g_1g_2\in E(G)$ and $h_1=h_2$, or $g_1=g_2$ and $h_1h_2\in E(H)$. All these products are associative and, with the exception of the lexicographic product, also commutative. For more on the products see~\cite{hik-2011}.

Let $G$ and $H$ be graphs and $*$  be one of the four graph products
under consideration. If $h\in V(H)$, then the set
$G^{h}=\{(g,h)\in V(G * H):\ g\in V(G)\}$ is a $G$-\emph{layer}. 
By abuse of notation we will also consider
$G^{h}$ as the corresponding induced subgraph. Clearly $G^{h}$ is
isomorphic to $G$ unless $*$ is the direct product in which case it
is an edgeless graph of order $|V(G)|$. For $g\in V(G)$, the
$H$-\emph{layer} $^g\!H$ is defined as $^g\!H
=\{(g,h)\in V(G * H):\ h\in V(H)\}$. We may again consider $^g\!H$
as an induced subgraph when appropriate. 

The rest of the paper is organized as follows. In the next section we first observe that  $\gamma_{gr}^t(G\times H) \geq \gamma_{gr}^t(G)\gamma_{gr}^t(H)$ holds for arbitrary graphs $G$ and $H$ and conjecture that actually the equality always holds. In the rest of the section we prove the conjecture for several special cases. In Section~\ref{sec:lexico} we consider the lexicographic product and express $\ggrt(G\circ H)$ in terms of related invariant of the factors $G$ and $H$. As a consequence, formulas for the Grundy total domination number of several special lexicographic products are obtained. In Section~\ref{sec:strong} lower bounds on $\ggrt(G \boxtimes H)$ are proved, while upper bounds are obtained for strong products of paths and cycles. The Cartesian product seems to be the most demanding with respect to the Grundy total domination number, a typical situation when domination problems are investigated on graph products, cf.~\cite{br-2012}. In Section~\ref{sec:cartesian} we then give upper and lower bounds for Cartesian products of paths and cycles. In the concluding section we briefly discuss two related invariants, L-Grundy domination number and Z-Grundy domination number and observe that some of the results derived in this paper extend to these two invariants.  

\section{Direct product}
\label{sec:direct}

We start this section with the following general bound on the Grundy total domination number of the direct product of graphs. 

\begin{lemma}\label{lem:greaterorequal}
If $G$ and $H$ are graphs, then $\gamma_{gr}^t(G\times H) \geq \gamma_{gr}^t(G)\gamma_{gr}^t(H)$.
\end{lemma}

\proof
Let $S=(x_1,\ldots,x_k)$ be a Grundy total dominating sequence in $G$, let $S'=(y_1,\ldots,y_{\ell})$ be a Grundy total dominating sequence in $H$, and let $x_i$ and $y_i$ totally footprint $x_i'$ and $y_i'$, respectively. Then $(x_1,y_1),(x_1,y_2),\ldots, (x_1,y_{\ell}),$ $(x_2,y_1),\ldots,(x_k,y_{\ell})$ is a total dominating sequence in $G\times H$ since $(x_i,y_j)$ totally footprints $(x_i',y_j')$. In fact, $(x_i',y_j') \in N(x_i,y_j)$, since $x_i' \in N_G(x_i)$ and $y_j' \in N_H(y_j)$, and if $(x_i',y_j') \in N((x_m,y_n))$ for some $m,n$, then $m\ge i$ and $n\ge j$. 
\qed

We conjecture that the inequality in Lemma~\ref{lem:greaterorequal} is in fact  equality.

\begin{conj}\label{conj:tdirect}
If $G$ and $H$ are graphs, then $\gamma_{gr}^t(G\times H) = \gamma_{gr}^t(G)\gamma_{gr}^t(H)$.
\end{conj}

In \cite{bhr-2016} a correspondence between total dominating sequences in graphs and edge covering sequences in hypergraphs was established. In the latter we are given sets (hyperedges) $E_1,\ldots, E_k$ that cover the ground set $V$ and we are looking for a maximum length $\rho_{gr}$ of a sequence $E_{i_1},\ldots, E_{i_{\ell}}$ such that $E_{i_m}\not \subseteq \bigcup_{n<m}E_{i_n}$. This problem clearly generalizes the problem of total dominating sequences and dominating sequences in graphs when sets $E_1,\ldots, E_k$ are the open or the closed neighborhoods of vertices, respectively.

For sets $X_1, X_2$ of hyperedges covering $V_1,V_2$, respectively, consider the product of $H_1 = (V_1, X_1)$, $H_2 = (V_2, X_2)$ as $H_1\times H_2$ with hyperedges $ \{E_1 \times E_2\mid E_1\in X_1, E_2 \in X_2 \}$ covering $V_1\times V_2$.  
Seemingly stronger conjecture in this setting is that $\rho_{gr}(H_1\times H_2) = \rho_{gr}(H_1)\rho_{gr}(H_2)$. To see that the conjectures are in fact equivalent, for each of the given $H_1,H_2$ construct the bipartite (incidence) graph $B_i$ whose vertices are $V_i \cup X_i$ and two vertices $u,v$  are adjacent if $v\in V_i$, $u \in X_i$ and $v\in u$. It is not hard to see that since $N((v_1,v_2)) = N(v_1)\times N(v_2)$, $\gamma_{gr}^t(B_1\times B_2)$ includes the questions of calculating $\rho_{gr}$ on four independent product structures, one of them being $H_1\times H_2$. If $\gamma_{gr}^t(B_1\times B_2)= \gamma_{gr}^t(B_1) \gamma_{gr}^t(B_2)$, then the conjecture holds for all the four subproblems. In particular, it holds that $\rho_{gr}(H_1\times H_2) = \rho_{gr}(H_1)\rho_{gr}(H_2)$.

In \cite{bgmrr-2014} a similar conjecture was posed stating that $\gamma_{gr}(G\boxtimes H) = \gamma_{gr}(G)\gamma_{gr}(H)$ and a construction connecting $\gamma_{gr}$ and $\rho_{gr}$ was introduced. Using the construction from the paper it can  be shown that this conjecture is in fact equivalent to both of the above conjectures. Moreover, examining the arguments used above, it suffices to prove Conjecture~\ref{conj:tdirect} by only considering that both factors are bipartite graphs. In fact, the conjecture holds for all pairs of bipartite graphs on at most 10 vertices, which was checked using a computer.

\begin{lemma}\label{lem:patition}
Let $E_1, \ldots, E_k$ be subsets of the edge set $E(G)$ of a graph $G$ such that $E_1\cup \dots \cup E_k=E(G)$. Let $G_1, \ldots, G_k$ be the isolate-free graphs with edge sets $E_1, \ldots, E_k$, respectively. Then $\gamma_{gr}^t(G)\leq \gamma_{gr}^t(G_1)+ \cdots + \gamma_{gr}^t(G_k)$.
\end{lemma}

\proof
Let $S=(x_1,\ldots,x_{\ell})$ be a Grundy total dominating sequence in $G$. For each $x_i\in S$ choose $y_i$ such that $x_i$ totally footprints $y_i$. Each pair $x_i,y_i$ induces an edge, thus it is in (at least) one of the graphs $G_1, \ldots, G_k$. If it is in more than one of them, choose one arbitrarily. We claim that for each $1\leq j \leq k$, the subsequence $S_j$ of those vertices $x_i\in \widehat{S}$ that are (chosen) in $G_j$, is a legal open neighborhood sequence in $G_j$. Indeed, when $x_i$ is chosen in $G_j$, it totally footprints $y_i$, for if $y_i$ was totally dominated in $G_j$ by $x_{r}, r<i$, then $y_i$ would be totally dominated also in $G$ by $x_{r}$, a contradiction. Hence, $\ggrt(G)=|\widehat{S}|=|\widehat{S_1}|+\cdots+|\widehat{S_k}|\le \ggrt(G_1)\cdots+\ggrt(G_k)$.
\qed

Let $bc(G)$ be the smallest size of a covering of edges of $G$ with complete bipartite graphs.

\begin{cor}
If $G$ is a graph, then $\gamma_{gr}^t(G)\leq 2bc(G)$.
\end{cor}

\proof
The result follows from Lemma~\ref{lem:patition} by using 
the fact that the Grundy total domination number of complete bipartite graphs is $2$.
\qed

\begin{lemma}\label{lem:H_v}
Let $v\in V(H)$, and let $S$ be a Grundy total dominating sequence in $G\times H$. Then $|\widehat{S}\cap G^v|\leq \gamma_{gr}^t(G)$.
\end{lemma}

\proof
Note that the subsequence of vertices in $S$ that lie in $G^v$ forms a legal open neighborhood sequence of $G\times H$, and consequently, its projection to $G$ is a legal open neighborhood sequence of $G$. Thus, $|\widehat{S}\cap G^v|\leq \gamma_{gr}^t(G)$.
\qed

It is known and easy to see that for any vertices $v_1,v_2$ with $N(v_1)=N(v_2)$ in $G$, we have $\gamma_{gr}^t(G) = \gamma_{gr}^t (G-v_2)$.

\begin{lemma}\label{lem:timesbipartite}
If $G$ is a graph, then  
$\gamma_{gr}^t(G\times K_{k_1,k_2})=2\gamma_{gr}^t(G)$. 
\end{lemma}

\proof
Note that any two vertices $u$ and $v$ in $K_{k_1,k_2}$ that are in the same bipartition set have the same open neighborhoods in $K_{k_1,k_2}$. Moreover, as $N_{G\times H}(u,v) = N_G(u) \times N_H(v)$, we infer that for any $g\in V(G)$, the vertices $(g,u)$ and $(g,v)$ have the same open neighborhoods in $G\times H$. By applying several times the observation before the lemma, we derive that $\gamma_{gr}^t(G\times K_{k_1,k_2})=\gamma_{gr}^t(G \times K_2)$. By Lemma \ref{lem:greaterorequal}, $\gamma_{gr}^t(G \times K_2)\geq 2\gamma_{gr}^t(G)$, and by Lemma~\ref{lem:H_v} the equality holds, that is, $\gamma_{gr}^t(G\times K_{k_1,k_2})=\gamma_{gr}^t(G\times K_2)=2\gamma_{gr}^t(G)$.
\qed

\begin{thm}\label{thm:equal}
If $G$ is a graph for which $\gamma_{gr}^t(G) = 2bc(G)$, then $\gamma_{gr}^t(G\times H) = \gamma_{gr}^t(G) \gamma_{gr}^t(H)$.
\end{thm}

\proof
Let $E_1, \ldots, E_k$ be the partition of the edge set of $G$ into complete bipartite graphs, such that $k=bc(G)$. Then $E_1\times E(H), \ldots, E_k\times E(H)$ is a partition of the edge set of $G\times H$ such that the corresponding graphs $G_1,\ldots,G_k$ are isomorphic to the direct product of complete bipartite graphs with $H$. By Lemma~\ref{lem:timesbipartite}, $\gamma_{gr}^t(G_i)=2\gamma_{gr}^t(H)$, thus by Lemma~\ref{lem:patition}, $\gamma_{gr}^t(G\times H) \leq 2bc(G) \gamma_{gr}^t(H)= \gamma_{gr}^t(G) \gamma_{gr}^t(H)$. By Lemma~\ref{lem:greaterorequal}, the equality holds.
\qed

\begin{cor}\label{cor:trees}
If $T$ is a tree, then $\gamma_{gr}^t(T\times H) = \gamma_{gr}^t(T) \gamma_{gr}^t(H)$.
\end{cor}

\proof
By Theorem \ref{thm:equal}, it suffices to show that $\gamma_{gr}^t(T) = 2 bc(T)$. By the result in \cite{bknt2017}, $\gamma_{gr}^t(T) = 2\beta(T)$, where $\beta(T)$ is the vertex cover number of $T$. Since the vertex cover can be interpreted as covering edges with stars, and stars are the only complete bipartite graphs in trees, it holds $\beta(T)=bc(T)$.
\qed

To find more graphs for which Theorem \ref{thm:equal} applies, we first consider the following lemma, which is a straightforward consequence of definitions.

\begin{lemma}\label{lem:twobipartite}
The direct product $K_{k_1,k_2}\times K_{\ell_1,\ell_2}$ of two complete bipartite graphs is isomorphic to the disjoint union $K_{k_1\ell_2, k_2\ell_2} + K_{k_1\ell_2,k_2\ell_1}$.
\end{lemma}

\begin{lemma}
If  $G_1,G_2$ are such that $\gamma_{gr}^t(G_i) = 2bc(G_i)$ for $i\in \{1,2\}$, then $\gamma_{gr}^t(G_1\times G_2) = 2bc(G_1 \times G_2)$.
\end{lemma}

\proof
On one hand, $ 2bc(G_1 \times G_2) \geq \gamma_{gr}^t(G_1\times G_2) = \gamma_{gr}^t(G_1)\gamma_{gr}^t(G_2)= 2bc(G_1)2bc(G_2)$. On the other hand, if $E_1, \ldots, E_{k_1}$ is a minimal partition of $G_1$ and $F_1, \ldots, F_{k_2}$ is a minimal partition of $G_2$, then by Lemma~\ref{lem:twobipartite}, each $E_i\times F_j$ is a disjoint union of two complete bipartite graphs. Thus $G_1 \times G_2$ has a partition of edges into $2bc(G_1)bc(G_2)$ complete bipartite graphs. This proves that 
the first inequality is in fact equality.
\qed

Another way to obtain graphs for which Theorem \ref{thm:equal} holds, is to take a graph $G$ with $\gamma_{gr}^t(G) = 2bc(G)$, copy a vertex in it and connect it to the neighbors of the original vertex. Then $\gamma_{gr}^t(G\times H) = \gamma_{gr}^t(G) \gamma_{gr}^t(H)$ for every graph $H$. For example, the latter holds if $G$ is a complete graph with a loop in each vertex.
Moreover, from results in~\cite{bgmrr-2014} it follows that the same   statement also holds if $G$ is a caterpillar with a loop in each vertex.

Notice that $\gamma^t_{gr}(P_\ell)=\ell, \gamma^t_{gr}(C_\ell)=\ell-2$ if $\ell$ is even, and $\gamma^t_{gr}(P_\ell)=\ell-1, \gamma^t_{gr}(C_\ell)=\ell-1$ otherwise. Hence, by applying Corollary~\ref{cor:trees}, we derive the following result. 

\begin{cor} If $k,\ell > 1$, then
\begin{displaymath}
\ggrt(P_k \times P_\ell)= \left\{ \begin{array}{l l}
{k}\cdot \ell, & \textrm{$k,\ell$ even}\\
{k}\cdot (\ell-1), & \textrm{$k$ even, $\ell$ odd}\\
(k-1)\cdot (\ell-1), & \textrm{$k,\ell$ odd}\\
\end{array}
\right.
\end{displaymath}
If $k> 1$, $\ell>2$, then
\begin{displaymath}
\ggrt(P_k \times C_\ell)= \left\{ \begin{array}{l l}
{k}\cdot (\ell-2), & \textrm{$k,\ell$ even}\\
{k}\cdot (\ell-1), & \textrm{$k$ even, $\ell$ odd}\\
(k-1)\cdot (\ell-2), & \textrm{$k$ odd, $\ell$ even}\\
(k-1)\cdot (\ell-1), & \textrm{$k,\ell$ odd}\\
\end{array}
\right.
\end{displaymath}
\end{cor}

\begin{thm}\label{thm:directcycles}
$\gamma^t_{gr}(C_{n_1}\times C_{n_2})=\gamma^t_{gr}(C_{n_1})\gamma^t_{gr}(C_{n_2})$.
\end{thm}

\proof
Notice that with a permutation of vertices, open neighborhoods of $C_\ell$ can be presented as $\{0,1\},\{1,2\},\{2,3\},\ldots,\{\ell-1,0\}$ in the case $\ell$ is odd and as a disjoint union of two such presentations in the case $\ell$ is even. For the sake of convenience, assume that both $n_1$ and $n_2$ are odd, we will deal with the other cases at the end of the proof. In this case the open neighborhoods of the product $C_{n_1}\times C_{n_2}$ are in the above presentation quadruples of the form $\{(0, 0), (1,0), (0,1), (1,1)\},\{(1, 0), (2,0), (1,1), (2,1)\},\ldots,\{(n_1-1, n_2-1), (0,n_2-1), (n_1-1,0), (0,0)\}$. We will call a set of vertices $\{ (i, 0), (i, 1), \ldots, (i, n_2 - 1)\}$ a row and similarly  a set of vertices $\{ (0,i), (1,i), \ldots, (n_1 - 1, i)\}$ a column.

Let $S=(v_1,\ldots,v_n)$ be an optimal total dominating sequence in  $C_{n_1}\times C_{n_2}$ of length $n$. For each $i\in [n] = \{1,\ldots,n\}$ let $k_i= |\bigcup_{j=1}^i N(v_j)| - |\bigcup_{j=1}^{i-1} N(v_j)| - 1$. By definition of the total dominating sequence, $k_i \geq 0$.
To prove that $\gamma^t_{gr}(C_{n_1}\times C_{n_2}) \leq \gamma^t_{gr}(C_{n_1})\gamma^t_{gr}(C_{n_2}) = n_1n_2 - n_1 - n_2 + 1$ holds, we need to prove that $\sum_{j=1}^{n} k_j \geq n_1 + n_2 - 1$. By Lemma \ref{lem:greaterorequal}, this suffices to prove the statement.

For the convenience denote $D_i = \bigcup_{j=1}^i N(v_j)$ for each $i \in [n]$. Let $v_i\in \{v_{2},\ldots, v_{n}\}$ be such that it totally dominates a vertex which is in a row or a column such that no vertex in this row or a column is in $D_{i-1}$. We will say that this row or column is newly dominated by $v_i$.

If $v_i$  newly dominates one column and not any row, it totally dominates at least two vertices not in $D_{i-1}$, therefore $k_{i}\geq 1$. The same holds if $v_i$ newly dominates one row and not any column.
 If $v_i$ newly dominates one column and newly dominates one row, it totally dominates at least three vertices not in $D_{i-1}$, therefore $k_{i}\geq 2$. Similarly, if $v_i$ newly dominates two columns (or rows) and one or zero rows (or columns), it totally dominates at least four vertices not in $D_{i-1}$, therefore $k_{i} = 3$. In all of the above cases when a row or a column is newly dominated, the corresponding $k_{i}$ contributes a positive value for each row and column to the final sum.
  
What remains is to analyze the situation where there exists $v_i$ that newly dominates two columns and two rows. Such vertex dominates four new vertices, hence $k_{i} = 3$. Vertex $v_1$ is certainly such, in the case that it is the only one, we already have that $\sum_{j=1}^{n} k_j \geq n_1 + n_2 - 1$ since the number of columns and rows is $n_1 + n_2$ and besides one row or column that is newly dominated by $v_1$ all of them contribute to the sum. To grasp the complexity of the case when more of the vertices that newly dominate two columns and two rows exist, we need additional definition.

For each $D_j$, let $G_j$ be a graph on vertices $D_j$ such that two vertices in $D_j$ are adjacent if they lie in the same column or the same row in $C_{n_1}\times C_{n_2}$. Notice that for $n \geq j_1 > j_2 \geq 1$, graph $G_{j_2}$ is an induced subgraph of $G_{j_1}$. Let $c_j$ be the number of connected components of $G_j$. If $v_i$, $i>1$, newly dominates two columns and two rows, then $c_i = c_{i-1} + 1$. Now let $v_{i'}$ be such that $c_{i'} < c_{i'-1}$, i.e. the number of connected components in $G_{i'}$ is less than in $G_{i'-1}$. There are a couple of cases to be analyzed here.

By definition of the graphs, $c_{i'-1} - c_{i'}$ can be at most 3 since vertices of $N(v_{i'})$ lie in two columns and two rows, thus they can join at most four components into one. First, suppose $c_{i'-1} - c_{i'}=3$. In this case each of the four components has a vertex in exactly one distinct row or column that vertices of $N(v_{i'})$ lie in. In particular, this implies that $v_{i'}$ does not newly dominate any row or column. Moreover, no vertex in $N(v_{i'})$ is in $D_{i'-1}$ since otherwise the component of such a vertex would lie in a column and a row that vertices of $N(v_{i'})$ lie in, hence there could not be four components in $G_{i'-1}$ that are joined into one in $G_{i'}$. Hence $v_{i'}$ must totally dominate four vertices not in $D_{i'-1}$ and $k_{i'} = 3$. 



In the case that $c_{i'-1} - c_{i'} = 2$ it is not hard to see that $k_{i'} \geq 2$, moreover it is important to notice that in this case it is possible that $v_{i'}$ also newly dominates one row or one column, but in this case $k_{i'} = 3$. Finally if $c_{i'-1} - c_{i'} = 1$ then similarly $k_{i'} \geq 1$, if $v_{i'}$ also newly dominates one column or one row then $k_{i'} \geq 2$, and if $v_{i'}$ also newly dominates two columns, two rows or a column and a row then $k_{i'} = 3$. In all of the above cases joining of connected components corresponds to adequate increase of $k_{i'}$.

We conclude the proof with the following. We have showed that when a column or a row is newly dominated, this is shown in an increase of the corresponding $k_i$, except in the case when there is $v_i$ that newly dominates two columns and two rows. In that case $k_i=3$ and, moreover, the number of connected components in $G_i$ is greater than in $G_{i-1}$, i.e. $c_{i} - c_{i-1} = 1$. Since $c_n = 1$ and $k_{i}$ is also adequately  increased if $c_{i-1} - c_{i} > 0$, it follows that $\sum_{j=1}^{n} k_j \geq n_1 + n_2 - 1$. This proves the assertion of the theorem when $n_1, n_2$ are odd.

To conclude the proof we need to cover the cases when at least one of the $n_1, n_2$ is even.
If, say, $n_1$ is odd and $n_2$ is even, then 
the open neighborhoods of the product $C_{n_1}\times C_{n_2}$ can be presented as a disjoint union of two copies of quadruples of the form $\{(0, 0), (1,0), (0,1), (1,1)\},\{(1, 0), (2,0), (1,1), (2,1)\},\ldots,\{(n_1-1, \frac{n_2}{2}-1), (0,\frac{n_2}{2}-1), (n_1-1,0), (0,0)\}$. Hence by similar arguments as above, 
$$\gamma^t_{gr}(C_{n_1}\times C_{n_2})\le 2(n_1-1)\left(\frac{n_2}{2}-1\right)=(n_1-1)(n_2-2)=\gamma^t_{gr}(C_{n_1})\gamma^t_{gr}(C_{n_2}).$$ 
Moreover, if $n_1$ and $n_2$ are even, we have four disjoint copies of quadruples of the form $\{(0, 0), (1,0), (0,1), (1,1)\},\{(1, 0), (2,0), (1,1), (2,1)\},\ldots,\{(\frac{n_1}{2}-1, \frac{n_2}{2}-1), (0,\frac{n_2}{2}-1), (\frac{n_1}{2}-1,0), (0,0)\}$ and get 
$$\gamma^t_{gr}(C_{n_1}\times C_{n_2})\le 4 \left(\frac{n_1}{2}-1 \right) \left(\frac{n_2}{2}-1 \right)=(n_1-2)(n_2-2)=\gamma^t_{gr}(C_{n_1})\gamma^t_{gr}(C_{n_2}).$$ 
\qed

Similar techniques (but more direct) can be used to show that Conjecture~\ref{conj:tdirect} holds if one of the factors is a cycle and one of them is a complete graph or both being complete graphs. This could possibly indicate why the conjecture holds for small graphs, since most of them could probably be partitioned into cycles, cliques, trees and complete bipartite graphs in the way of Lemma \ref{lem:patition}.

\section{Lexicographic product}
\label{sec:lexico}

As it turns out, the formula for $\ggrt(G\circ H)$ a bit surprisingly relies on (Grundy) dominating
sequences of $G$. The results are very similar to the formula for
$\ggr(G\circ H)$ which was obtained in~\cite{bbgkkptv2016}.

Given a dominating sequence $D=(d_1,\ldots,d_k)$ in a graph $G$, let
$a(D)$ denote the cardinality of the set of vertices $d_i$ from $D$,
which are not adjacent to any vertex from $\{d_1,\ldots,d_{i-1}\}$.

\begin{thm}
\label{thm:lex} For any graphs $G$ and $H$, where $H$ has no isolated
vertices:
 $$\ggrt(G\circ H)=
\max\{a(D)(\ggrt(H)-1)+|\widehat{D}|\,;\,D\textit{ is a dominating
sequence of $G$}\}.$$
\end{thm}

\proof Let $D=(d_1,\ldots,d_m)$
be a dominating sequence of $G$, and let $(d'_1,\ldots,d'_k)$ be a
Grundy total dominating sequence of $H$.
Then one can find a sequence $S$ in $G\circ H$ of length
$a(D)(\ggrt(H)-1)+|\widehat{D}|$ as follows. Let $S$ be the
sequence that corresponds to $D$, and those vertices $d_i\in
\widehat{D}$ which are not adjacent to any vertex from
$\{d_1,\ldots,d_{i-1}\}$ are repeated $\ggrt(H)$
times in a row, so that the corresponding subsequence is of the form
$((d_i,d'_1),\ldots,(d_i,d'_k))$. On the other hand, the vertices
$d_i\in \widehat{D}$ which are adjacent to some vertex from
$\{d_1,\ldots,d_{i-1}\}$ are projected only once from the vertices
of $S$, i.e. there is a unique vertex $(d_i,h)$ that belongs to $S$.
It is easy to see that in either case the vertices in $S$ are
legally chosen. In the first case this is true because no vertex of
$^{d_i}\!H$ is totally dominated yet at
the point when $(d_i,d'_1)$ is chosen, thus
$((d_i,d'_1),\ldots,(d_i,d'_k))$ is a legal open neighborhood
subsequence. In the second case this
is true because $D$ is a legal closed neighborhood sequence in $G$,
and so when $d_i$ is chosen, there exists another vertex $t\in
V(G)$, distinct from $d_i$, that $d_i$ footprints. Hence, when
$(d_i,h)$ is chosen in $S$, it totally footprints
vertices from $^{t}\!H$. Note that the length of $S$ is
$a(D)(\ggrt(H)-1)+|\widehat{D}|$. This implies that $\ggrt(G\circ
H)\ge \max\{a(D)(\ggrt(H)-1)+|\widehat{D}|\,;\,D\textit{ is a
dominating sequence of $G$}\}$.

For the converse, let $S$ be an arbitrary open neighborhood sequence
in $G\circ H$. Let
$S'=\big((x_1,y_1),\ldots,(x_n,y_n)\big)$ be the subsequence of $S$,
where $(x,y)\in S'$ if and only if $(x,y)$ is the first vertex in
$S$ that belongs to $^{x}\!H$. We claim that the corresponding
sequence of the first coordinates $T=(x_1,\ldots,x_n)$ is a legal
closed neighborhood sequence in $G$. Firstly, if $x_i$ is not
adjacent to any of $\{x_1,\ldots,x_{i-1}\}$, then clearly $x_i$
footprints itself. Otherwise, if $x_i$ is adjacent to $x_j$, where
$j\in\{1,\ldots,{i-1}\}$, then there exists a vertex $(x_j,h)\in S$
that appears in $S$ before any vertex from $^{x_i}\!H$. Thus when
$(x_i,y_i)$ is added to $S$, all vertices from $^{x_i}\!H$ are
already totally dominated. Hence, since
$(x_i,y_i)$ must totally footprint some vertex,
it can only be a vertex in $^{x}\!H$, where $x\in N_G(x_i)$. We
infer that $x_i$ footprints $x$ with respect to $T$. Thus $T$ is a
legal closed neighborhood sequence.

Let $A(T)$ be the set of all vertices $x_i$ in $T$, such that $x_i$
is not adjacent to any vertex from $\{x_1,\ldots,x_{i-1}\}$. (Note
that $|A(T)|=a(T)$ by definition.) Two cases for a vertex $x_j\in
\widehat{T}$ appear, which give different bounds on the number of
vertices in $^{x_j}\!H\cap\widehat{S}$. If $x_j\notin A(T)$, then
$|^{x_j}\!H\cap\widehat{S}|=1$, because $(x_j,y_j)$ totally
dominates all neighboring layers, and $^{x_j}\!H$ is
already totally dominated before $(x_j,y_j)$ is added
to $S$. On the other hand, if $x_j\in A(T)$, then clearly
$|^{x_j}\!H\cap\widehat{S}|\le \ggrt(G)$. We infer that
$|\widehat{S}|\le (|\widehat{T}|-a(T))+a(T)\ggrt(H)$, where $T$ is
a closed neighborhood sequence. \qed

Since any independent set of vertices yields a legal closed
neighborhood sequence, we infer the following
\begin{cor}
For any graphs $G$ and $H$ with no isolated vertices:
$$\ggrt(G\circ H)\ge {\alpha}(G)\ggrt(H).$$
\end{cor}

\begin{figure}[ht!]
\centering
\begin{tikzpicture}[scale=.9,style=thick]
\def\vr{2.5pt}


\coordinate (p1) at (-1,0); \coordinate (p2) at (0,0); \coordinate
(p3) at (1,1); \coordinate (p4) at (2,1.4); \coordinate (p5) at
(2,0.6); \coordinate (p6) at (1,-1); \coordinate (p7) at (2,-1.4);
\coordinate (p8) at (2,-0.6);

\draw (p1) -- (p2) -- (p3) -- (p4); \draw (p3) -- (p5); \draw (p2)
-- (p6) -- (p7); \draw (p6) -- (p8);


\draw(p1)[fill=black] circle(\vr); \draw(p2)[fill=white]
circle(\vr); \draw(p3)[fill=white] circle(\vr);
\draw(p4)[fill=black] circle(\vr); \draw(p5)[fill=black]
circle(\vr); \draw(p6)[fill=white] circle(\vr);
\draw(p7)[fill=black] circle(\vr); \draw(p8)[fill=black]
circle(\vr);

\node[below] at (p2){$v$}; \node[below] at (p1){$u$};

\end{tikzpicture}
\vskip -0.3 cm \caption{The tree $T_8$.} \label{fig:tree}
\end{figure}
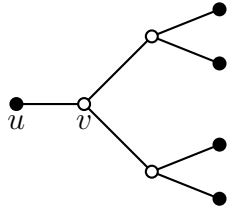

To see that the inequality in the corollary above is not always equality, let $T_8$
be the tree from Fig.~\ref{fig:tree}, and let $H$ be a graph with
no isolated vertices. Note that the filled vertices in the figure
present a maximum independent set of $T_8$, and so $\alpha(T_8)=5$.
Hence the bound in the corollary is $\alpha(T_8)\ggrt(H)=5\ggrt(H)$, 
but $\ggrt(T_8\circ H)\ge 5\ggrt(H)+1$. Indeed, let $S$ be the sequence
that starts in the layer $^u\!H$, by legally picking $\ggrt(H)$
vertices, then choosing a vertex from $^v\!H$, and then
choosing $\ggrt(H)$ vertices in each of the layers that
correspond to the remaining four leaves of $T_8$. This yields a legal
sequence of the desired length.

Following the same arguments as in~\cite[Corollary
10]{bbgkkptv2016}, and using the fact that $\ggrt(H)\ge 2$
for any graph $H$ with no isolated vertices, we
derive the following results.

\begin{cor}\label{cor:PnH} If $H$ is a graph with no isolated
vertices, then
\begin{displaymath}
\ggrt(P_k \circ H)= \left\{ \begin{array}{l l}
\frac{k}{2}\cdot \ggrt(H)+1, & \textrm{$k$ is even, $k\neq 2$}\\
\left\lceil \frac{k}{2} \right\rceil \cdot \ggrt(H), & \textrm{$k$ is odd.}\\
\end{array}
\right.
\end{displaymath}
\end{cor}

\begin{cor} Let $k,\ell > 2$.
If $\ell$ is even, then
\begin{displaymath}
\ggrt(P_k \circ P_\ell)= \left\{ \begin{array}{l l}
\frac{k}{2}\cdot \ell+1, & \textrm{$k$ is even}\\
\left\lceil \frac{k}{2} \right\rceil \cdot \ell, & \textrm{$k$ is odd.}\\
\end{array}
\right.
\end{displaymath}

If $\ell$ is odd, then
\begin{displaymath}
\ggrt(P_k \circ P_\ell)= \left\{ \begin{array}{l l}
\frac{k}{2}\cdot (\ell-1)+1, & \textrm{$k$ is even}\\
\left\lceil \frac{k}{2} \right\rceil \cdot (\ell-1), & \textrm{$k$ is odd.}\\
\end{array}
\right.
\end{displaymath}
\end{cor}

\begin{cor} Let $k,\ell > 2$.
If $\ell$ is even, then
\begin{displaymath}
\ggrt(P_k \circ C_\ell)= \left\{ \begin{array}{l l}
\frac{k}{2}\cdot (\ell-2)+1, & \textrm{$k$ is even}\\
\left\lceil \frac{k}{2} \right\rceil \cdot (\ell-2), & \textrm{$k$ is odd.}\\
\end{array}
\right.
\end{displaymath}

If $\ell$ is odd, then
\begin{displaymath}
\ggrt(P_k \circ C_\ell)= \left\{ \begin{array}{l l}
\frac{k}{2}\cdot (\ell-1)+1, & \textrm{$k$ is even}\\
\left\lceil \frac{k}{2} \right\rceil \cdot (\ell-1), & \textrm{$k$ is odd.}\\
\end{array}
\right.
\end{displaymath}
\end{cor}

Using Theorem~\ref{thm:lex} we also obtain the following result.

\begin{cor}
Let $H$ be a graph without isolated vertices, and $k\ge 3$. Then
\begin{displaymath}
\ggrt(C_k \circ H)= \left\{ \begin{array}{l l}
\frac{k}{2}\cdot \ggrt(H), & \textrm{$k$ is even}\\
\left\lfloor  \frac{k}{2} \right\rfloor \cdot \ggrt(H) +1, & \textrm{$k$ is odd.}\\
\end{array}
\right.
\end{displaymath}
\end{cor}

The value for Grundy total domination numbers of
$C_k\circ C_\ell$ easily follows from the above corollary.

\section{Strong product}
\label{sec:strong}

Given a (closed neighborhood) dominating sequence $D=(d_1,\ldots , d_k)$ in a graph $G$, let 
\begin{itemize}
\item $C(D)$ denote the set of vertices $d_i$ from $D$ that footprints itself and at least one of its neighbors and let $c(D)=|C(D)|$. 
\item $B(D)$ denote the set of vertices $d_i$ from $D$ that does not footprint itself and footprint at least one of its neighbors and let $b(D)=|B(D)|$.
\item $A(D)$ denote the set of vertices $d_i$ from $D$ that are not adjacent to any vertex from $\{d_1,\ldots , d_{i-1}\}$ and let $a(D)=|A(D)|$.
\end{itemize}

\begin{thm} 
If $G$ and $H$ are graphs, where $G$ has no isolated vertices, 
\begin{eqnarray*}
     \begin{aligned}
     &\ggrt(G \boxtimes H) \geq \max \big\{ (\ggr(G)+1)\cdot c(D)+\ggr(G)\cdot b(D)+(|\widehat{D}|-b(D)-c(D))\ggrt(G);\\
        & \hspace{37mm} D \textrm{ is a dominating sequence of $H$}\big\}.
     \end{aligned}
\end{eqnarray*}
\end{thm}

\proof
Let $D=(y_1,\ldots , y_k)$ be an arbitrary legal closed neighborhood sequence of $H$. Let $(x_1,\ldots , x_\ell)$ be a $\ggr$-sequence of $G$ and $(g_1,\ldots , g_j)$ a $\ggrt$-sequence of $G$. Let $x \in V(G)$ be a vertex that is footprinted by $x_\ell$ (note that $x_\ell$ can be chosen in such a way that $x\neq x_\ell$, since $G$ has no isolated vertices). Now we construct a legal total dominating sequence $S=S_k$ of $G \boxtimes H$. First let $S_0 = \emptyset$. For each $i \in [k]$ let $S_i=S_{i-1} \oplus ((x_1,y_i),(x_2,y_i),\ldots ,(x_{\ell-1},y_i),(x,y_i),(x_\ell,y_i))$ if $y_i \in C(D)$. If $y_i \in B(D)$, let $S_i=S_{i-1} \oplus ((x_1,y_i),(x_2,y_i),\ldots ,(x_{\ell-1},y_i), (x_\ell,y_i))$ otherwise $S_i=S_{i-1} \oplus ((g_1,y_i),(g_2,y_i),\ldots ,(g_j,y_i))$.  
\qed

By symmetry we get the following.

\begin{thm} 
If $G$ and $H$ are graphs, where $H$ has no isolated vertices, 
\begin{eqnarray*}
     \begin{aligned}
     &\ggrt(G \boxtimes H) \geq \max \big\{ (\ggr(H)+1)\cdot c(D)+\ggr(H)\cdot b(D)+(|\widehat{D}|-b(D)-c(D))\ggrt(H);\\
        & \hspace{37mm} D \textrm{ is a dominating sequence of $G$}\big\}.
     \end{aligned}
\end{eqnarray*}
\end{thm}

\begin{thm}
If $G$ and $H$ are graphs, then 
$$\ggrt(G \boxtimes H) \geq \max{\{\ggrt(G)\cdot a(D)+(|\widehat{D}|-a(D))\ggr(G); D \textrm{ is a dominating sequence of $H$}\}}. $$
\end{thm}

\proof
Let $D=(y_1,\ldots, y_k)$ be an arbitrary legal closed neighborhood sequence of $H$. Let $(x_1,\ldots, x_\ell)$ be a $\ggr$-sequence of $G$ and $(g_1,\ldots, g_j)$ a $\ggrt$-sequence of $G$. Now we construct a legal total dominating sequence $S=S_k$ of $G \boxtimes H$. Let $S_0 = \emptyset$. For each $i \in [k]$ let $S_i=S_{i-1} \oplus ((g_1,y_i),(g_2,y_i),\ldots ,(g_j,y_i))$ if $y_i \in A(D)$, otherwise $S_i=S_{i-1} \oplus ((x_1,y_i),(x_2,y_i),\ldots ,(x_\ell,y_i))$.  
\qed

By symmetry we also get the following.

\begin{thm}
If $G$ and $H$ are graphs, then 
$$\ggrt(G \boxtimes H) \geq \max{\{\ggrt(H)\cdot a(D)+(|\widehat{D}|-a(D))\ggr(H); D \textrm{ is a dominating sequence of $G$}\}}. $$
\end{thm}

For the strong product of paths and/or cycles the above theorems imply the following lower bounds.
\begin{cor} For any pair $k, \ell \ge 3$ of integers we have
 \begin{enumerate}
				\item[(i)]
				\begin{displaymath}
				\ggrt(P_k \boxtimes P_\ell) \geq \left\{ \begin{array}{l l}
					k\ell-k-\ell+2, & \textrm{$k,\ell$ are odd}\\
					k\ell-k-\ell+3, & \textrm{otherwise.}\\
				\end{array}
				\right.
				\end{displaymath}

				\item[(ii)]
				\begin{displaymath}
				\ggrt(P_k \boxtimes C_\ell) \geq \left\{ \begin{array}{l l}
					k\ell-2k-\ell+3, & \textrm{$k$ is odd, $\ell$ is even}\\
					k\ell-2k-\ell+4, & \textrm{otherwise.}\\
				\end{array}
				\right.
				\end{displaymath}

				\item[(iii)]
				\begin{displaymath}
				\ggrt(C_k \boxtimes C_\ell) \geq \left\{ \begin{array}{l l}
					k\ell-2k-2\ell+5, & \textrm{$k,\ell$ are even}\\
					k\ell-2k-2\ell+6, & \textrm{otherwise.}\\
				\end{array}
				\right.
				\end{displaymath}

    \end{enumerate}
\end{cor}

We round this section with upper bounds on the Grundy total dominations number of strong products of paths and/or cycles. 

\begin{thm}\label{upperStrong} For any pair $k, \ell \ge 3$ of integers we have
    \begin{enumerate}
				\item[(i)]
        $\ggrt (P_k \boxtimes P_\ell) \le k\ell-\min\{k,\ell\}$,
        \item[(ii)]
        $\ggrt (C_k \boxtimes C_\ell)\le k\ell-\min\{2k,2\ell\}$,        
				\item[(iii)]
				$\ggrt(P_k \boxtimes C_\ell) \le k\ell- \min\{2k,\ell\}$.
    \end{enumerate}
\end{thm}

\proof
Let $V(P_k) = [k]$, $E(P_k)=\{12,23,\ldots, (k-1)k\}$.
Note that if an open neighborhood sequence intersects every $P_k$-layer in less than $k$ vertices and every $P_\ell$-layer in less than $\ell$ vertices, then the length of the sequence is at most $k\ell-\max\{k,\ell\}$. 

\begin{figure}[ht!]
    \begin{center}
        \begin{tikzpicture}[scale=1.0,style=thick,x=1cm,y=1cm]
        \def\vr{2.5pt} 

    
    \draw (-7,0) -- (-5,0);
		\draw (-3,0)--(-2,0);
		\draw (2,0)--(4,0);
		\draw (6,0)--(7,0);
		\draw (0,0)--(0,2);
		\draw (0,4)--(0,5);

        \foreach \x  in {-7,-5,-2,0,2,4,7}
        {
        \filldraw [fill=black, draw=black,thick] (\x,0) circle (3pt);
        }
				 \foreach \x  in {-6,-3,3,6}
        {
        \filldraw [fill=white, draw=black,thick] (\x,0) circle (3pt);
        }
				
				\foreach \y in {1,2,4,5}
				\filldraw [fill=black, draw=black,thick] (0,\y) circle (3pt);

		\draw (-7,-0.3) node {$x_1$};
		\draw (-6,-0.3) node {$a_1$};
		\draw (-5,-0.3) node {$x_2$};
		\draw (-4,0) node {$\ldots$};
		\draw (-3,-0.3) node {$a_p$};
		\draw (-2,-0.3) node {$y_1$};
		\draw (-1,0) node {$\ldots$};
    \draw (1,0) node {$\ldots$};
		\draw (2,-0.3) node {$y_r$};
		\draw (3,-0.3) node {$a_{p+1}$};
		\draw (4,-0.3) node {$x_{p+1}$};
		\draw (5,0) node {$\ldots$};
		\draw (6,-0.3) node {$a_o$};
		\draw (7,-0.3) node {$x_o$};
		\draw (-8,0) node {$P_k^i$};
		\draw (0,3) node {$\vdots$};
		\draw (0,5.5) node {$^jP_\ell$};

        \end{tikzpicture}
    \end{center}
    \caption{Layers $P_k^i$ and $^{j}P_\ell$.}
    \label{proof}
\end{figure}
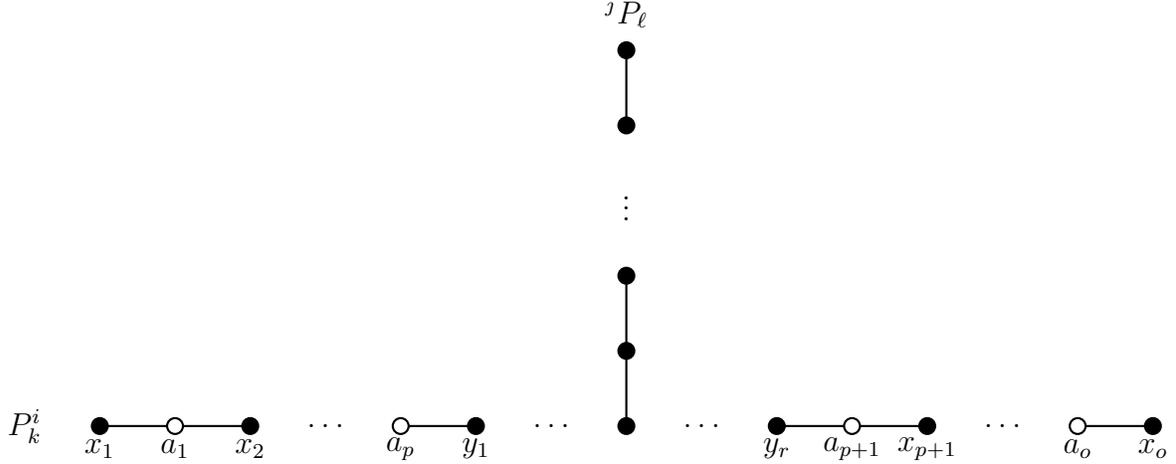

Consider now the smallest $m$ such that $D=(x_1,\dots,x_m)$ contains all vertices from a $P_k$-layer or all vertices from a $P_\ell$-layer. If $x_1,\dots,x_m$ would contain one complete $P_k$-layer and one complete $P_\ell$-layer, then we would have $N(x_m)\subseteq \cup_{i=1}^{m-1}N(x_i)$, a contradiction. Assume without loss of generality that the sequence $x_1,\dots,x_m$ contains all vertices of a $P_\ell$-layer $^jP_\ell$. Let $n_i$, $i\in [\ell]$, be the number of vertices from the $i^{\rm th}$ $P_k$-layer that are in $x_1,\dots,x_m$. Then at least $n_i+1$ vertices from $P_k^i$ are totally dominated by $\widehat{D} \cap P_k^i$ and hence footprinted by $D$, unless vertices from $\widehat{D} \cap P_k^i$ are of the form $\{1,3,\ldots , j-a,j-a+1,\ldots , j+b, j+b+2,j+b+4, \ldots , k-2,k\}$. If we use notation of Figure~\ref{proof}, $\widehat{D} \cap P_k^i=\{x_1,x_2,\ldots , x_o\} \cup \{y_1,\ldots , y_r\}$. In this case, if there exists $i' \in \{1,\ldots , o\}$ such that $x_{i'} \in P_k^* \cap \widehat{D}$ or $a_{i'} \in P_k^* \cap \widehat{D}$, for $* \in \{i-1,i+1\}$, then at least $n_i+1$ vertices from $P_k^i$ are footprinted by $D$. If for any $i' \in \{1,\ldots , o\}$, $a_{i'},x_{i'} \notin P_k^* \cap \widehat{D}$, then at least $n_i +2o$ vertices from $P_k^*$ are footprinted by $D$. As $o \geq 2$, $m$ vertices from $D$ footprint at least $\sum_{i=1}^\ell (n_i+1)=m+\ell$ vertices of $P_k \boxtimes P_\ell$. Therefore any legal open neighborhood sequence in $P_k \boxtimes P_\ell$ contains at most $k\ell-\ell$ vertices.

The proof for $C_k \boxtimes C_\ell$ proceeds in a similar way. First we note that if $D$ is a legal open neighborhood sequence, then $|C_k^i \cap \widehat{D}| \leq k-1$ and $|^jC_\ell \cap \widehat{D}| \leq \ell-1$ for any $i \in [\ell], j \in [k]$. 
Consider now the first moment $m$ when $D=(x_1,\dots,x_m)$ contains $k-1$ vertices from a $C_k$-layer or $\ell-1$ vertices from a $C_\ell$-layer. Assume without loss of generality that the sequence $x_1,\dots,x_m$ contains $\ell-1$ vertices of a $C_\ell$-layer $^jC_\ell$. Let $n_i$, $i\in [\ell]$, be the number of vertices from the $i^{\rm th}$ $C_k$-layer that are in $x_1,\dots,x_m$. Then we can prove in the similar way as in $(i)$, that $n_i+2$ vertices of each $C_k^i$ layer are footprinted by $D$.

The proof of $(iii)$ goes in a similar way, we have to consider two cases depending on either it happens first that legal open neighborhood sequence contains $k$ vertices from a $P_k$-layer or $\ell-1$ vertices from a $C_\ell$-layer. \qed

\section{Cartesian product}
\label{sec:cartesian}

\begin{thm}\label{upperCart} For any pair $k, \ell \ge 3$ of integers we have
    \begin{enumerate}
        \item[(i)]
        $\ggrt (C_k \cp C_\ell)\le k\ell-\min\{k,\ell\}$,
        \item[(ii)]
        $\ggrt (P_k \cp P_\ell) \le k\ell-\min\{\lfloor k/2 \rfloor,\lfloor \ell/2 \rfloor\}$,
        \item[(iii)]
        $\ggrt (P_k \cp C_\ell) \le k\ell- \min\{k,\lceil \ell/2 \rceil \}$.
    \end{enumerate}
\end{thm}

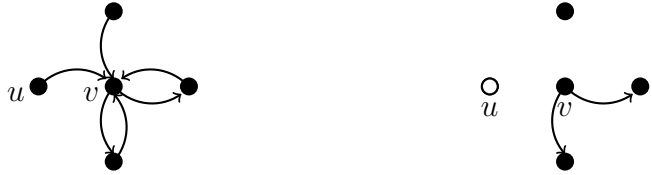
\begin{figure}[ht!]
    \begin{center}
        \begin{tikzpicture}[scale=1.0,style=thick,x=1cm,y=1cm]
        \def\vr{2.5pt} 

        \foreach \x  in {-4,-3,-2,2,3,4}
        {
					\filldraw [fill=black, draw=black,thick] (\x,0) circle (3pt);
        }

    \foreach \x  in {-3,3}
        {
        \filldraw [fill=black, draw=black,thick] (\x,-1) circle (3pt);
        }

				\foreach \x  in {-3,3}
        {
        \filldraw [fill=black, draw=black,thick] (\x,1) circle (3pt);
        }
				
				\filldraw [fill=white, draw=black,thick] (2,0) circle (3pt);

    \draw (-4.3,-0.1) node {$u$}; 
		\draw (-3.3,-0.1) node {$v$}; 
		\draw (3,-0.3) node {$v$};
		\draw (2,-0.3) node {$u$};
   
		\draw[->] (-4,0) to [bend left=40] (-3.1,0.1);
		\draw[->] (-2,0) to [bend right=40] (-2.9,0.1);
		\draw[->] (-3,1) to [bend right=40] (-3,0.1);
		\draw[->] (-3,-1) to [bend right=40] (-3,-0.1);
		\draw[->] (-3,0) to [bend right=40] (-3,-0.9);
		\draw[->] (-3,0) to [bend right=40] (-2.1,-0.1);
		\draw[->] (3,0) to [bend right=40] (3,-0.9);
		\draw[->] (3,0) to [bend right=40] (3.9,-0.1);

        \end{tikzpicture}
    \end{center}
    \caption{The left part of the figure presents the situation before $u$ is added to the sequence, the right part presents the situation after $u$ is added to the sequence.}
    \label{fig:Removal}
\end{figure}

\proof (i) Consider the directed graph $G$ obtained from $C_k\cp C_\ell$ by replacing every edge with two opposite arcs. We interpret the arc from $u$ to $v$ in $G$ as a certificate that $u$ can be added to an open neighborhood sequence $x_1,\dots,x_m$ as $v\notin \cup_{i=1}^mN(x_i)$. As the open neighborhood sequence grows, we remove those arcs of $G$ for which the above statement is no longer valid (see Figure~\ref{fig:Removal}). The number of arcs in $G$ is $4k\ell$ out of which $2k\ell$ are in $C_k$-layers and the same number in $C_\ell$-layers. We claim that whenever we add a vertex to an open neighborhood sequence, we have to remove at least two arcs in $C_k$-layers and at least two arcs in $C_\ell$-layers. Indeed,
since $u$ is added to the open neighborhood sequence because it dominates a neighbor $v$ that was not dominated before, then the arcs from all four neighbors of $v$ to $v$ should be removed.

Note that if an open neighborhood sequence intersects every $C_k$-layer in less than $k$ vertices and every $C_\ell$-layer in less than $\ell$ vertices, then the length of the sequence is at most $k\ell-\max\{k,\ell\}$. Moreover, if $x_1,\dots,x_m$ would intersect one complete $C_k$-layer and one complete $C_\ell$-layer with $x_m$ being the intersection of these layers, then since $k,\ell \ge 3$ we would have $N(x_m)\subseteq \cup_{i=1}^{m-1}N(x_i)$, a contradiction. 

Consider now the first moment $m$ when $x_1,\dots,x_m$ contains all vertices from a $C_k$-layer or all vertices from a $C_\ell$-layer (but not both!). Assume without loss of generality that the former occurs and the sequence $x_1,\dots,x_m$ contains all vertices of the $j^{\rm th}$ $C_k$-layer. Let $n_i$, $i\in [k]$, be the number of vertices from the $i^{\rm th}$ $C_\ell$-layer that are in $x_1,\dots,x_m$. Note that all $2n_i$ arcs outgoing from these $n_i$ vertices have been removed. Moreover, if $(v_i,u_s),\dots, (v_i,u_p)$ is a maximal sequence of consecutive vertices such that all of them are contained in $\{x_1,\dots,x_m\}$ and either $s\neq p$ or $s=j=p$ holds, then also the arcs incoming to $(v_i,u_s)$ and $(v_i,u_p)$ have been deleted.
As at least one vertex of this $C_\ell$-layer is not in $x_1,\dots,x_m$, the number of arcs removed in the layer is at least $2n_i+2$. Summing up we obtain that the number of arcs removed in $C_\ell$-layers is at least $\sum_{i=1}^{k}\left(2n_i+2\right) =  2(m+k)$. According to our observation above, in every later step at least two arcs of $C_\ell$-layers will be removed, therefore the length of the complete open neighborhood sequence is at most $$m+\frac{2k\ell-2(m+k)}{2}=k\ell-k$$ and we are done with the proof of (i).

\vspace{3mm}

(ii) and (iii) We modify the proof above a little. We introduce the weighted directed  graphs $G'$ and $G''$ obtained from $P_k\cp P_\ell$ and $P_k\cp C_\ell$ respectively, by replacing every edge with two opposite arcs. All arcs have weight 1 except $((v_2,u_j),(v_1,u_j))$ and $((v_{k-1},u_j),(v_k,u_j))$ for all $j\in [\ell]$ in both $G'$ and $G''$, and in $G'$ also $((v_i,u_2),(v_i,u_1))$ and $((v_i,u_{\ell-1}),(v_i,u_\ell))$ for all $i\in [k]$. These arcs have weight 2. In this way, the total weight of the arcs in $P_k$-layers and $P_\ell$-layers (resp.\ $C_\ell$-layers) is $2k\ell$ both in $G'$ and $G''$. It remains also true that when adding a vertex to an open neighborhood sequence then both the arcs removed  in $P_k$-layers and $P_\ell$-layers (resp.\ $C_\ell$-layers) have total weight at least 2.

    The rest of the proof is very similar to that of part (i). Let $m$ be the first moment when a complete $P_k$-layer or $P_\ell$-layer (resp.\ $C_\ell$-layer) of $P_k \cp P_\ell$ or $P_k\cp C_\ell$ belongs to the open neighborhood sequence $x_1,\ldots, x_m$. Note that this time it is possible that the sequence $x_1,\ldots, x_m$ contains both a $P_k$-layer and a $P_\ell$-layer (resp.\ $C_\ell$-layer).

\medskip\noindent
\textsc{Case I:} $x_1,\ldots, x_m$ contains only a $P_k$-layer or a $P_\ell$-layer (resp.\ $C_\ell$-layer). \\
If in $P_k\cp C_\ell$ the sequence $x_1,\ldots, x_m$ contains a $P_k$-layer, then with the same proof as in (i) we obtain $\ggrt (P_k \cp C_\ell)\le k\ell-k$. So we can assume that the sequence $x_1,\ldots, x_m$ contains a $C_\ell$-layer and we will prove $\ggrt (P_k \cp C_\ell)\le k\ell-\lceil \ell/2 \rceil$. (The case of $P_\ell$-layer for $P_k\cp P_\ell$ is analogous.)

    If the number of vertices in the $j^{th}$ $P_k$-layer is $m_j$, then we have $1\le m_j < k$ for all $1 \le j \le \ell$. Therefore the total weight of arcs removed in $P_\ell$-layers from $G'$ or $G''$ is at least $2m_j+1$ with equality if and only if the $m_j$ vertices form a subpath of the path containing one of its endpoints. Summing up for all $P_k$-layers we obtain that the total weight of removed arcs in $P_k$-layers is at least $2m+\ell$. Therefore the length of the complete open neighborhood sequence is at most $$m+\frac{2k\ell-(2m+\ell)}{2}.$$

\medskip\noindent
\textsc{Case II:} $x_1,\ldots, x_m$ contains both a $P_k$-layer and a $P_\ell$-layer (resp.\ $C_\ell$-layer). \\
    We consider the open neighborhood sequence $x_1,\ldots, x_{m-1}$. As $x_1,\ldots, x_m$ contains both a $P_k$-layer and a $P_\ell$-layer (resp.\ $C_\ell$-layer), we know that every $P_k$-layer and $P_\ell$-layer (resp.\ $C_\ell$-layer) contains at least one vertex from $x_1,\ldots, x_{m-1}$ and one vertex not from $x_1,\ldots, x_{m-1}$. Therefore the argument can be finished as in the previous case.
\qed

\begin{thm}
For any pair $k,\ell \geq 3$ we have
\begin{enumerate}
\item $\ggrt(P_k \Box P_\ell) \geq k\ell-\min{\{k,\ell\}}$,
\item $\ggrt(P_k\Box C_\ell) \geq k\ell - \min{\{2k,\ell\}}$,
\item $\ggrt(C_k \Box C_\ell) \geq k\ell- \min{\{2k,2\ell\}}$.
\end{enumerate}
\end{thm}

\proof
Let the vertex set of $X \cp Y$ be the following:

$$V(X \cp Y):=\{(a,b) : 0 \le a \le k-1,\ 0 \le b \le \ell-1\},$$

\noindent
and let the edge set be:
$$E(X \cp Y):=\{((a,b),(c,d)) \in V(X \cp Y) \times V(X \cp Y):$$
$$(a=c \pm 1 \ (a \equiv c \pm 1 \textrm{ (mod $k$) if $X=C_k$) and } b=d ) \textrm{ or }$$
$$( a=c \textrm{ and } b= d \pm 1 \ (b \equiv d \pm 1 \textrm{ (mod  $\ell$) if $Y=C_\ell$))} \}.$$

\vspace{2mm}

\noindent
We also define the lexicographic and antilexicographic ordering on $V(X \cp Y)$. For any $(a,b),(c,d) \in V(X \cp Y)$ let

$$(a,b) \prec_{lex} (c,d) \textrm{ if and only if } a < c \textrm{ or } (a=c \textrm{ and }b < d),$$
$$(a,b) \prec_{alex} (c,d) \textrm{ if and only if } b < d \textrm{ or } (b=d \textrm{ and }a < c).$$

\vspace{2mm}

Let $v_1,\ldots, v_{k\ell}$ be the lexicographic ordering of vertices of $X \cp Y$, with $|V(X)|=k, |V(Y)|=\ell$. Furthermore, let $u_1,u_2,\ldots, u_{k\ell}$ be the antilexicographic ordering of vertices of $X \cp Y$.
Then $(v_1,v_2,\ldots , v_{(k-1)\ell})$ and $(u_1,u_2,\ldots , u_{(\ell-1)k})$ are legal total dominating sequences of $P_k \cp P_\ell$. Thus $\ggrt(P_k \cp P_\ell) \geq k\ell-\min{\{k,\ell\}}$.

On the other hand $(v_1,\ldots , v_{(k-1)\ell})$ and $(u_1,\ldots , u_{(\ell-2)k})$ are legal total dominating sequences of $P_k \cp C_\ell$. Thus $\ggrt(P_k \cp C_\ell) \geq k\ell-\min{\{2k,\ell\}}$.

Finally, $(v_1,\ldots , v_{(k-2)\ell})$ and $(u_1,\ldots , u_{(\ell-2)k})$ are legal total dominating sequences of $C_k \cp C_\ell$. Thus $\ggrt(C_k \cp C_\ell) \geq k\ell-\min{\{2k,2\ell\}}$.
\qed

We note that the following is known \cite{D2014,IMA2010} for the Grundy total domination number of the Cartesian product of two paths of the same length: for $k\ge 1$ we have $$\ggrt(P_k \cp P_k)=k^{2}-k.$$ 

On the other hand, for an odd $k$ it holds that $C_k \times C_k \cong C_k \cp C_k$, hence by Theorem \ref{thm:directcycles} it holds $\ggrt(C_k \times C_k) = k^2 - 2k + 1$.
\section{Concluding remarks}
 
In~\cite{bbgkkptv2017} two additional versions of dominating sequences were introduced and studied. When the definition~\eqref{eq:defGrundyT} is modified to read as
\begin{equation}
\label{eq:defGrundyZ}
N(v_i) \setminus \bigcup_{j=1}^{i-1}N[v_j]\not=\emptyset\,,
\end{equation}
we get the so-called Z-sequences, and when
\begin{equation}
\label{eq:defGrundyL}
N[v_i] \setminus \bigcup_{j=1}^{i-1}N(v_j)\not=\emptyset\,,
\end{equation}
the so-called L-sequences are defined provided that vertices in such a sequence are distinct. The corresponding invariants obtained from the longest possible (Z- or L-) sequences are denoted by $\ggrz(G)$ and $\ggrl(G)$, respectively. As it turns out, these two invariants have natural counterparts also in the zero-forcing and minimum-rank world, see~\cite{L-2017}. 

With a little effort, by using almost the same proofs as in this paper, one can prove some of the results for the remaining two invariants.

\begin{itemize}
\item Version of Theorem~\ref{thm:lex} holds also for the $L$-Grundy domination number of lexicographic products:
\begin{thm} For any graphs $G$ and $H$ with no isolated
vertices:
 $$\ggrl(G\circ H)=
\max\{a(D)(\ggrl(H)-1)+|\widehat{D}|\,;\,D\textit{ is a dominating
sequence of $G$}\}.$$
\end{thm}
Clearly, one can also derive exact values of this invariant in lexicographic products of paths and/or cycles. 

\item Lemma~\ref{lem:greaterorequal} (and its proof) holds also if Grundy total dominaton number is replaced with the $Z$-Grundy domination number, i.e., $$\gamma_{gr}^Z(G\times H)\geq \gamma_{gr}^Z(G)\gamma_{gr}^Z(H)$$ holds for any graphs $G$ and $H$. However, for $\gamma_{gr}^z$ the left-hand side can be strictly greater, as demonstrated by $K_3\times K_3$, where $\gamma_{gr}^z(K_3)=1$ while $\gamma_{gr}^z(K_3\times K_3)=4$.

\item Also, Lemma~\ref{lem:patition} can be proved in the setting of any of the four Grundy domination invariants. 
\begin{lemma}
Let $E_1, \ldots, E_k$ be the subsets of the edge set of a graph $G$. Let $G_1, \ldots, G_k$ be graphs on $E_1, \ldots, E_k$, respectively. Then $\gamma_{gr}^x(G)\leq \gamma_{gr}^x(G_1)+ \cdots + \gamma_{gr}^x(G_k)$ for $x \in \{Z,t,L,\emptyset\}$.
\end{lemma}
\end{itemize}

\subsection*{Acknowledgement}

The research of Cs.\ Bujt\'as, B.\ Patk\'os, Zs.\ Tuza and M.\ Vizer was supported by the National Research, Development and Innovation Office -- NKFIH, grant SNN 116095.

The authors acknowledge the project (Combinatorial Problems with an Emphasis on Games,
N1-0043) was financially supported by the Slovenian Research Agency. The authors acknowledge the financial support from the Slovenian Research Agency (research core funding No.\ P1-0297). G. Ko\v smrlj was also supported by Slovenian Research Agency under the grant P1-0294 and Slovenian Public Agency for Entrepreneurship, Internationalization, Foreign Investments and Technology under the grant KKIPP-99/2017.



\begin{thebibliography}{99}

\bibitem{Aazami2008} 
  A.~Aazami, 
  Hardness results and approximation algorithms for some problems on graphs, 
  PhD thesis, University of Waterloo, 2008, http://hdl.handle.net/10012/4147.

\bibitem{AIM} 
  AIM Minimum Rank-Special Graphs Work Group,
  Zero-forcing sets and the minimum rank of graphs,
  \textit{Linear Algebra Appl.} \textbf{428} (2008) 1628--1648.

\bibitem{ajps-2016} 
  T.~Ansill, B.~Jacob, J.~Penzellna, D.~Saavedra,
  Failed skew zero forcing on a graph, 
  \textit{Linear Algebra Appl.} \textbf{509} (2016) 40--63. 

\bibitem{BBF2013} 
  F.~Barioli, W.~Barrett, S.~M.~Fallat, H.~T.~Hall, L.~Hogben, B.~Shader, P.~van~den~Driessche, H.~van~der~Holst,
  Parameters related to tree-width, zero forcing, and maximum nullity of a graph, \textit{J.\ Graph Theory} \textbf{72} (2013) 146--177.

\bibitem{BFF-2016+}
  K.~F.~Benson, D.~Ferrero, M.~Flagg, V.~Furst, L.~Hogben, V.~Vasilevskak, B.~Wissman, 
  Power domination and zero forcing,
   arXiv:1510.02421v4 [math.CO] (22 Feb 2017). 

\bibitem{bbgkkptv2016} 
  B.~Bre{\v{s}}ar, Cs.~Bujt\'as, T.~Gologranc, S.~Klav\v{z}ar, G.~Ko{\v{s}}mrlj, B.~Patk\'os, Zs.~Tuza, M.~Vizer,   
  Dominating sequences in grid-like and toroidal graphs,
  \textit{Electron.\ J.\ Combin.} \textbf{23}(4) (2016) P4.34.

\bibitem{bbgkkptv2017} 
  B.~Bre{\v{s}}ar, Cs.~Bujt\'as, T.~Gologranc, S.~Klav\v{z}ar, G.~Ko{\v{s}}mrlj, B.~Patk\'os, Zs.~Tuza, M.~Vizer,
  Grundy dominating sequences and zero forcing sets,
  \textit{Discrete Optim.} \textbf{26} (2017) 66--77. 
  
\bibitem{br-2012}
  B.~Bre\v{s}ar, P.~Dorbec, W.~Goddard, B.~Hartnell, M.~A.~Henning, S.~Klav\v{z}ar, D.~F.~Rall,
  Vizing's conjecture: a survey and recent results,
\textit{J.\ Graph Theory} \textbf{69} (2012) 46--76.

\bibitem{bknt2017} 
  B.~Bre{\v{s}}ar, T.~Kos, G.~Nasini, P.~Torres. 
  Total dominating sequences in trees, split graphs, and under modular decomposition.
  \textit{Discrete Optim.}, to appear.


\bibitem{bgmrr-2014}
  B.~Bre{\v{s}}ar, T.~Gologranc, M.~Milani\v c, D.~F.~Rall, R.~Rizzi,
  Dominating sequences in graphs,
  \textit{Discrete Math.} \textbf{336} (2014) 22--36.

\bibitem{bhr-2016}
  B.~Bre{\v{s}}ar, M.~A.~Henning, D.~F.~Rall, 
  Total dominating sequences in graphs,
  \textit{Discrete Math.} \textbf{339} (2016) 1665--1676.

\bibitem{D2014} 
  L.~M.~DeAlba,
  Some results on minimum skew zero forcing sets, and skew zero forcing number,
  arXiv preprint arXiv:1404.1618. (2014)

\bibitem{hik-2011}
  R.~Hammack, W.~Imrich, S.~Klav\v{z}ar,
  Handbook of Product Graphs, Second Edition,
  CRC Press, Boca Raton, FL, 2011.


    
\bibitem{heye-2013}
  M.~A.~Henning, A.~Yeo, 
  Total Domination in Graphs, 
  Springer, 2013.


\bibitem{IMA2010} 
  IMA-ISU research group on minimum rank (M. Allison, E. Bodine, L. M. DeAlba, J. Debnath, L. DeLoss, C. Garnett, J. Grout, L. Hogben, B. Im, H. Kim, R. Nair, O. Pryporova, K. Savage, B. Shader, and A. WangsnessWehe),
  Minimum rank of skew-symmetric matrices described by a graph,
  \textit{Linear Algebra Appl.} \textbf{432} (2010) 2457--2472.

\bibitem{K2015} 
  N.~F.~Kingsley,
  Skew propagation time,
  PhD thesis, Iowa State University, 2015, http://lib.dr.iastate.edu/etd/14804/

\bibitem{L-2017}
  J.~C.-H.~Lin, 
  Zero forcing number, Grundy domination number, and their variants, 
  arXiv:1706.00798 (June 6, 2017).
  
\bibitem{masu-2016} 
  S.~Mallik, B.~L.~Shader, 
  On graphs of minimum skew rank 4, 
  \textit{Linear Multilinear Algebra} \textbf{64} (2016) 279--289. 

\end{thebibliography}
\end{document}